\documentclass{article}

\usepackage[utf8]{inputenc}
\usepackage[english]{babel}

\usepackage{graphicx}
\usepackage{amssymb}
\usepackage{amsmath}
\usepackage{tikz-cd}
\usepackage{enumerate}

\usepackage{hyperref}

\newtheorem{theorem}{Theorem}[section]

\newtheorem{proposition}[theorem]{Proposition}

\newcommand{\C}{\mathcal{C}}
\newcommand{\U}{\mathcal{U}}

\begin{document}

\title{Deligne categories as limits in rank and characteristic}
\author{Nate Harman}

\maketitle

\begin{abstract}
We give new interpretations of the Deligne categories $\underline{Rep}(GL_t)$ and $\underline{Rep}(S_t)$ (and their abelian envelopes) over $\mathbb{C}$ in terms of modular representations of general linear and symmetric groups of large rank in large characteristic. In particular we make sense of the sentence ``$\underline{Rep}(S_n)$ is the limit of $Rep(S_{p+n})$ over $\overline{\mathbb{F}}_p$ as $p$ goes to infinity". We then give examples of how to pass results between these different settings.

\end{abstract}

\begin{section}*{Introduction}

Deligne constructed rigid symmetric Karoubian tensor categories $\underline{Rep}(GL_t)$ and $\underline{Rep}(S_t)$ interpolating the categories of representations of general linear groups (symmetric groups respectively), which are universal for having an object (resp. Frobenius algebra) of dimension $t$, where $t$ is an element of an arbitrary ring $k$ \cite{Deligne3}.  If $k$ is a field of characteristic zero these categories are semisimple and hence abelian unless $t$ is an integer (resp. nonnegative integer).   If $t$ is an integer (resp. nonnegative integer) Deligne constructed abelian tensor categories $\underline{Rep}(GL_t)^{ab}$ and $\underline{Rep}(S_t)^{ab}$ containing $\underline{Rep}(GL_t)$ and $\underline{Rep}(S_t)$, which were later shown to be universal for this property (in \cite{EHS} for $GL_t$, \cite{CO2} for $S_t$).

In this paper, we give a new interpretation of these tensor categories (both the Karoubian and abelian versions) in characteristic zero as limits of categories of representations of in positive characteristic where we let both the rank of the group and the characteristic of the field grow.  This then allows us to pass results from representation theory in positive characteristic to Deligne's categories and vice-versa.

In the first section we introduce a fairly general construction due to Deligne which uses an ultrafilter to construct interesting new categories out of sequences of related categories coming from representation theory, often capturing their asymptotic behavior. We give two warm up examples of this construction in action, and then apply it to our setting to relate these Deligne categories to modular representations.

Then in the second section we give examples of how to extract concrete results out of this abstract set-up. We give explicit formulas for irreducible modular representations of symmetric groups, a characterization of simple commutative algebras in the Deligne categories, and a connection between representations of general linear groups in positive characteristic and representations of supergroups in characteristic zero.

The purpose of this note is to introduce this interpretation and give a flavor for how it can be used, as such the examples are chosen to be simple and direct applications of this. We hope to apply these ideas to carry out more sophisticated calculations in subsequent work.

\end{section}

\begin{section}*{Acknowledgments}

Thanks to Pavel Etingof for many helpful conversations, and to Vera Serganova for her comments and feedback. This work was partially supported by the National Science Foundation Graduate Research Fellowship under Grant No. 1122374.

\end{section}

\begin{section}{Ultrafilters and limits of categories}

We will assume some basic familiarity with ultrafilters, including taking ultraproducts of rings, and the ultralimit of a sequence in a compact topological space.  For a reference on some general results about ultrafilters and ultraproducts see \cite{Schoutens}.

Let $\mathcal{U}$ be a non-principal ultrafilter on the set of natural numbers (which we will henceforth refer to just as an ultrafilter).  Given a family of categories\footnote{For technical reasons one should restrict to the class of small categories for this construction, but it is easy to then include essentially small categories, which includes all of the representation theoretic categories we care about. We will henceforth ignore such issues.} $\C_1, \C_2, \dots$  we may define their ultraproduct $\widehat{\C}_\U$ as follows:

Objects of $\widehat{\C}_\U$ are germs of families of object $X_n \in \C_n$ relative to $\U$. That is, sequences of objects $X_n \in \C_n$ where two sequences are the same if they agree on a subset of $\mathbb{N}$ belonging to $\U$.  Morphisms are given by germs of morphisms, defined in the same manner.  Composition is given by choosing representative sequences and performing the composition termwise, and one can show the germ of the composition is independent of the choice of representation.

If we impose extra structure on the categories $\C_i$ then we often get similar structure on $\widehat{\C}_\U$.  For example if the categories $\C_i$ are all abelian then so is $\widehat{\C}_\U$, as one may take representative families and define kernels and cokernels ``pointwise" in each $\C_n$. Similarly if each $\C_n$ has a monoidal structure, then so does $\widehat{\C}_\U$. 

In general, the category $\widehat{\C}_\U$ will be much too large to work with it (the ultraproduct construction does \emph{not} preserve most finiteness conditions).  So often we will restrict to a more manageable subcategory $\C_\U$ inside $\widehat{\C}_\U$; the choice of this subcategory will vary a bit with the situation.  We will demonstrate this with two examples:

\medskip

\noindent \textbf{Example 1 (Lefschetz principle for finite groups):} Let $G$ be a fixed finite group and take $\C_n$ to be the category of finite dimensional representations of $G$ over $\overline{\mathbb{F}}_{p_n}$ where $p_n$ denotes the $n$th prime number. Then $\widehat{\C}_\U$ is naturally a rigid symmetric tensor category over $\Pi_\U \overline{\mathbb{F}}_{p_n} \cong \mathbb{C}$  (fix such an isomorphism, we know one exists by Steinitz' theorem).  A first guess might be that $\widehat{\C}_\U$ is equivalent to the category of finite dimensional complex representations of $G$, but this is false since, in particular, $\widehat{\C}_\U$ has objects of infinite length (for example, consider the sequence $X_n \in \C_n$ where $X_n$ is a direct sum of $n$ copies of the trivial representation).

We can remedy this by letting $\C_\U \subset \widehat{\C}_\U$ be the subcategory generated by the regular representation (i.e. the object in $\widehat{\C}_\U$ represented by the sequence of regular representations) under taking direct sums and summands; then indeed it's easy to verify that $\C_\U$ is equivalent to the category of finite dimensional complex representations of $G$ (under our identification $\Pi_\U \overline{\mathbb{F}}_{p_n} \cong \mathbb{C}$).

We see dependence on $\U$ is very minimal; indeed for a different choice of ultrafilter $\U'$ we have that $\C_\U$ and $\C_{\U'}$ are equivalent as symmetric monoidal categories.  One needs to be careful at the level of sequences of objects and morphisms, since the same underlying sequence may correspond to very different germs under the two ultrafilters. Moreover, the complex structure depends heavily on ultrafilter as well as the choice of identification $\Pi_\U \overline{\mathbb{F}}_{p_n} \cong \mathbb{C}$.  Nevertheless, this independence on $\U$ does allow one to conclude that any first order statement (expressible in the language of symmetric monoidal categories) about representations of $G$ over $\mathbb{C}$ holds over algebraically closed fields of large enough characteristic.

\medskip 

\noindent \textbf{Example 2 (Deligne's letter to Ostrik):}  Fix a prime number $p$ and let $\C_n$ be the category of finite dimensional representations of the symmetric group $S_n$ (Deligne also considered general linear, orthogonal, and symplectic groups) over $\mathbb{F}_p$. It follows that $\widehat{\C}_\U$ is a rigid symmetric tensor category over $\Pi_\U \mathbb{F}_p \cong \mathbb{F}_p$, however it has infinite dimensional hom-spaces and objects of infinite length and is hence not pre-tannakian. 

 We can trim this down to a pre-tannakian category $Rep(S_\U)$ by taking the abelian tensor subcategory generated by the object $X$ represented by the sequence of standard $n$-dimensional permutation representations of $S_n$.  These categories were constructed as the first examples of pre-tannakian categories in positive characteristic with an object of super-exponential growth. 
 
 Deligne conjectured that these categories $\C_\U$ should depend on the choice of $\U$ very minimally. More precisely, he conjectured they would only depend on the ultralimit $t \in \mathbb{Z}_p$ of the sequence $1,2,3,\dots$ of natural numbers viewed as a sequence in $\mathbb{Z}_p$.  See more about this construction, and the resolution of this conjecture in the case of symmetric groups for $p \ge 5$ in \cite{Harman}.
 
 \medskip
 
 \noindent \textbf{Remark:} The weak dependence on the choice of ultrafilter suggests that we should be able to make some of the results constructive, and indeed for these examples one can.  In the first example one can pass results between characteristic $p$ for $p \nmid |G|$ and characteristic zero by lifting to Witt vectors.  In the second example, one can remove the need for an ultrafilter in Deligne's construction by using the categorical stabilization results for symmetric groups in \cite{Harman} and taking an appropriate direct limit. 

\medskip

The first example can be thought of as taking a ``limit in characteristic", and the second example can be thought of as taking a ``limit in rank" of these categories of representations.  Informally, the main idea of this note is that for symmetric and general linear groups it makes sense to take a ``simultaneous limit in characteristic and rank", and the resulting limit categories are just Deligne's categories $\underline{Rep}(GL_t)$ and $\underline{Rep}(S_t)$ over $\mathbb{C}$, or their corresponding abelian envelopes. The following main theorem makes this precise:

\begin{theorem}
Let $t_1, t_2, t_3, \dots$ be a sequence of positive integers tending to infinity, and let $p_1, p_2, p_3, \dots$ be a sequence of prime numbers tending to infinity. Let $t \in \mathbb{C}$ denote the image of the sequence $t_1 \in \mathbb{F}_{p_1} , t_2 \in \mathbb{F}_{p_2}, \dots$ under a fixed identification $\Pi_\U \overline{\mathbb{F}}_{p_n} \cong \mathbb{C}$.  Take $\C_n$ to be either $Rep(GL(t_n, \overline{\mathbb{F}}_{p_n}))$ or $Rep(S_{t_n}, \overline{\mathbb{F}}_{p_n})$, let $X_n \in \C_n$ denote the defining $t_n$-dimensional representation of $GL(t_n, \overline{\mathbb{F}}_{p_n})$ ($S_{t_n}$ respectively), and let $X \in \widehat{\C}_\U$ be the object corresponding to the sequence $X_1, X_2, X_3, \dots$. Then we have:

\begin{enumerate}[(a)]
\item The full subcategory $\C_\U \subset \widehat{\C}_\U$ generated by $X$ under the operations of taking duals, tensor products, direct sums, and direct summands is equivalent (as symmetric tensor categories over $\mathbb{C}$) to the (Karoubian) Deligne category $\underline{Rep}(GL_t)$ (resp. $\underline{Rep}(S_t)$).

\item The full subcategory $\C_\U^{ab} \subset \widehat{\C}_\U$ generated by $X$ under the operations of taking duals, tensor products, direct sums, and subquotients is equivalent to the abelian envelope $\underline{Rep}(GL_t)^{ab}$ (resp. $\underline{Rep}(S_t)^{ab}$) of the corresponding Deligne category.

\end{enumerate}

\end{theorem}

\noindent \textbf{Proof:}  The ultrafilter construction ensures that $\widehat{\C}_\U$ is a $\Pi_\U \overline{\mathbb{F}}_{p_n} \cong \mathbb{C}$-linear rigid symmetric tensor category, and the categorical dimension of $X$ is the ultralimit of the categorical dimensions of the $X_n$, which is $t$ by definition (and $X$ is a Frobenius algebra in the case of symmetric groups).  Hence the universal property of $\underline{Rep}(GL_t)$  (resp. $\underline{Rep}(S_t)$) we have a symmetric tensor functor from $\underline{Rep}(GL_t)$ (resp. $\underline{Rep}(S_t)$) into $\widehat{\C}_\U$ sending the defining representation to $X$, which lands inside $\C_\U$ by definition.

Symmetric tensor functors out of $\underline{Rep}(GL_t)$ are characterized in \cite{EHS} (and in \cite{CO2} for  $\underline{Rep}(S_t)$), and it is seen that the map we get from the universal property is faithful and extends to an exact tensor functor from $\underline{Rep}(GL_t)^{ab}$ (resp. $\underline{Rep}(S_t)^{ab}$) if and only if the object (resp. Frobenius algebra) $X$ is not killed by any Schur functor.  In our case, the condition that the sequence $t_n$ tends to infinity ensures that any Schur functor can only kill finitely many of the $X_n$ and hence doesn't kill $X$.

So all that remains is to verify that these functors are full.  For part (a), since all objects are direct summands of (direct sums of) objects of the form $X^{\otimes r}\otimes X^{*\otimes s}$ we can reduce the statement that the functor is an equivalence of Karoubian categories to the statement that $End(X^{\otimes r}\otimes X^{*\otimes s})$ is the appropriate walled Brauer algebra (resp. partition algebra) over $\mathbb{C}$ with parameter $t$. Schur-Weyl duality (which holds over any infinite field) tells us that the corresponding endomorphism algebra in $\mathcal{C}_n$ is the appropriate walled Brauer algebra (resp. partition algebra) over $\overline{\mathbb{F}}_{p_n}$ with parameter $t_n$, and it is then an easy exercise to show that the ultraproduct is then the correct algebra over $\mathbb{C}$. 

Finally, for part (b) we observe that the functor we get from $\underline{Rep}(GL_t)^{ab}$ (resp. $\underline{Rep}(S_t)^{ab}$) into $\C_\U^{ab}$ is a functor between highest weight categories over the same poset, which by part (a) is an equivalence when we restrict to the subcategories of tilting modules. Hence it follows that this is an equivalence of categories. $\square$

\medskip

\noindent \textbf{Remark:} Instead of taking a sequence of primes $p_i$ we may just apply the ultrafilter construction to the corresponding categories of representations over an algebraically closed field of characteristic zero.  This case is addressed by Deligne in the introduction of \cite{Deligne3}, it gives an interpretation of these Deligne categories at transcendental $t$ (one cannot get algebraic $t$ in this manner), and indeed Deligne uses this as motivation for his more algebraic description of these categories.  Really the main new observation here is that if we let the characteristic vary along with the rank then we can view these Deligne categories at algebraic (including integer) values of $t$ as limits in this sense.

\end{section}

\begin{section}{Applications, predictions, and open problems}

We now have some equivalences of categories which depend on a choice of ultrafilter as well an abstract isomorphism between $\mathbb{C}$ and a field of characteristic zero defined via the ultrafilter.  At first this may seem too abstract and non-constructive to be able to extract anything meaningful, so we'd now like to give a few examples of how if we choose the input to our main theorem carefully we can extract some concrete results and predictions in representation theory.

\begin{subsection}{Low degree modular characters of symmetric groups}

A basic question in the representation theory of symmetric groups is to characterize the irreducible representations of low dimension.  More precisely one could ask for some fixed $k$, what are all of the irreducible representations of $S_n$ of dimension at most $n^k$?

In characteristic zero this question was first considered by Rasala in \cite{Rasala} where he showed that for fixed $k$ and large $n$ the only irreducible $S_n$ representations (up to twisting by the sign character) of dimension at most $n^k$ are those Specht modules $S^{\lambda}$ corresponding to partitions $\lambda$ with $n-\lambda_1 \le k$.  Moreover he used the hook length formula to give explicit formulas for these dimensions as polynomials in $n$.

In \cite{James} James considered the problem in characteristic $p$, where he showed that for fixed $k$ and large $n$ the only irreducible representations of $S_n$ over a field of characteristic $p$ (up to twisting by the sign character) are those irreducible modules $D^\lambda$ corresponding to $p$-regular partitions $\lambda$ with $n-\lambda_1 \le k$.  However he was not able in general to give explicit formulas for their dimensions.

In \cite{BK} Brundan and Kleshchev gave explicit dimension formulas for these irreducible representations the case when $k=4$, $p \ge 5$, and $n$ is sufficiently large.  Interestingly, these formulas are independent of $p$ in a certain sense. For example they show that for $p \ge 5$ and $n$ large enough:

\begin{displaymath}
   \text{dim }D^{(n-2,2)} = \left\{
     \begin{array}{lr}
       (n^2-5n+2)/2 & : n \equiv 2 \mod p\\
       (n^2-3n-2)/2 & : n \equiv 1 \mod p\\
       (n^2-3n)/2 & : \text{otherwise} \hspace{.72cm}
     \end{array}
   \right.
\end{displaymath} 

Really what they calculated was how to express irreducibles in terms of Specht modules at the level of Grothendieck groups; these dimension formulas then follow from the polynomial dimension formulas for Specht modules given by Rasala.  Going through their proof, the above example really reflects the stronger statement that inside the Grothendieck group of $Rep(S_n)$ we have (for $p\ge5$ and $n$ large) that:

\begin{displaymath}
   [D^{(n-2,2)}] = \left\{
     \begin{array}{lr}
       [S^{(n-2,2)}]-[S^{(n-1,1)}]& : n \equiv 2 {\mod p}\\
      {[S^{(n-2,2)}]-[S^{(n)}]} & : n \equiv 1 {\mod p}\\
       {[S^{(n-2,2)}]} & :  \text{otherwise} \hspace{.72cm}
     \end{array}
   \right.
\end{displaymath} 

We will now generalize this for all $k$ and sufficiently large characteristics $p$. In particular we will express irreducibles $D^\lambda$ with $n-\lambda_1 \le k$ in terms of Specht modules inside the Grothendieck groups of the categories of representations. Combining this with dimension formulas for Specht modules gives exact dimension formulas for all irreducible representations of $S_n$ of dimension at most $n^k$ in characteristic $p$, provided $p$ and $n$ are large enough relative to $k$.

Fix $k$ along with another non-negative integer $m$.  For an integer $n$ and a partition $\lambda$ of $n$ satisfying $n-\lambda_1 \le k$.  Consider the set $A_{\lambda,m}$ of partitions consisting of $\lambda$ along with all partitions obtainable by first removing a $(n-m)$-rim hook from $\lambda$, and then re-adding one in a way that the result is higher than $\lambda$ in the dominance order.  If $n$ is sufficiently large (relative to $k$) then $A_\lambda$ inherits a linear order from the dominance order. We denote its elements in decreasing order as $\lambda^{(0)}, \lambda^{(1)}, \dots, \lambda^{(d)}$, in particular $\lambda^{(d)} = \lambda$.

\begin{proposition}
Suppose $p \gg k$, and $n > p$ is congruent to $m$ modulo $p$. Then the Specht module $S^\lambda$ over a field of characteristic $p$ is reducible iff $|A_{\lambda,m}| \ge 2$, and in that case in the Grothendieck group of $Rep(S_n)$ we have:

$$[D^\lambda] = [S^{\lambda^{(d)}}]-[S^{\lambda^{(d-1)}}]+ \dots +(-1)^d[S^{\lambda^{(0)}}]$$

\end{proposition}

\noindent \textbf{Proof:}  The periodicity results for decomposition matrices in \cite{Harman} ensure that for $p$ large enough it is sufficient to verify the claim for when $n = p+m$.  Then we note that if $t_i = m + p_i$ for some sequence of primes $p_i$ tending to infinity that $lim_\U t_i = m \in \mathbb{C} \cong \Pi_\U \overline{\mathbb{F}}_{p_n}$ independently of the choices of ultrafilter and field isomorphism.

We can then conclude the expression of $D^{(p+m-|\mu|, \mu)}$ in terms of Specht modules in $Rep(S_{p+m})$ over $\overline{\mathbb{F}}_p$ agrees with the corresponding expression for $D^{(m-|\mu|, \mu)}$ in the Deligne category $\underline{Rep}(S_m)^{ab}$ for all but finitely many primes $p$, since if we had infinitely many ``bad" primes  $p_i$ we could apply the main theorem to the sequence $t_i =m+p_i$ to pass the result to $Rep(S_m)^{ab}$.

Finally, Comes and Ostrik's description in \cite{CO} of blocks in $Rep(S_m)^{ab}$ tells us explicitly how to decompose a Specht module into irreducibles in the Grothendieck group.  Our formula for an irreducible in terms of Specht modules follows immediately from this by induction. $\square$

\medskip

The ultrafilter machinery is well suited to give asymptotic results, but not very good at telling us when the asymptotic behavior will set in.  One may be interested knowing how large $p$ and $n$ need to be for the above proposition to hold, and this machinery gives no explicit bounds.  We have the following conjecture in this direction:

\medskip

\noindent \textbf{Conjecture:} The conclusion of the above proposition holds whenever $p > k$ and $n > 4k$.

\medskip

We suspect this can be shown by careful analysis of the modular branching rule as in the proof in \cite{BK}, or by using the periodicity results in \cite{Harman} to reduce the problem to blocks of small defect where James' conjecture on decomposition numbers is known to hold and then analyzing the LLT algorithm.  Since this is not the focus of this paper we will not pursue this direction any further at this time.

\end{subsection}

\begin{subsection}{Representation theory in complex rank}

Etingof has introduced (in \cite{Etingof}) a program he calls ``representation theory in complex rank" in which he takes representation theoretic objects which are built out of symmetric or general linear group representations with additional structure, and tries to understand corresponding objects in $\underline{Rep}(S_t)$ or $\underline{Rep}(GL_t)$.

Often one can deduce the behavior of these categories at transcendental values of $t$ from the behavior at large integer values, however there is a recurring technical hurdle in this theory when one tries to characterize things at algebraic values of $t$.  We believe that the ideas presented here can help overcome this hurdle, and will demonstrate this with a simple example.

\medskip

In \cite{Luke} Sciarappa investigated simple commutative algebras inside the Deligne category $\underline{Rep}(S_t)$.  He showed that the only simple commutative algebras in $\underline{Rep}(S_t)$ for transcendental $t$ are those interpolating algebras of functions on cosets of $H \times S_{n-k}$ inside $S_n$ for some a fixed subgroup $H \subset S_k$. However his techniques were not able to eliminate the possibility of ``exotic"  simple commutative algebras at algebraic values of $t$.  We can now complete his characterization:

\begin{proposition}
Sciarappa's classification of simple commutative algebras inside $\underline{Rep}(S_t)$ holds for all values of $t\in \mathbb{C}$, as well as in $\underline{Rep}(S_t)^{ab}$ for integer values of $t$. 
\end{proposition}

\noindent \textbf{Proof:} As a step in his proof, Sciarappa showed that for large enough $n$ these algebras of functions on cosets of $H \times S_{n-k}$ are exactly those simple commutative algebras in $Rep(S_n)$ with underlying object is in the abelian subcategory of $Rep(S_n)$ generated by the defining representation and its first $k$ tensor powers.  We first note that this result holds over any algebraically closed ground field, regardless of the characteristic.

Suppose there was an exotic simple commutative algebra $A$ in $\underline{Rep}(S_t)$ for some algebraic number $t$, let $q(x)$ denote the minimal polynomial of $t$ (over the integers) with positive leading coefficient. Now we claim that we can find a sequence of integers $t_1, t_2, \dots$ and primes $p_1, p_2, \dots$ (both tending to infinity) such that $q(t_n) \equiv 0 \mod p_n$ for all $i$. To see this, it is enough to show that there are infinitely many prime numbers dividing some number of the form $q(n)$ for $n \in \mathbb{N}$. 

This then follows from an easy density argument; for large $N$ there are roughly $N^{1/deg(q)}$ positive integers of the form $q(n)$ with $q(n) < N$, whereas for a fixed collection of $m$ primes there are on the order of $\log(N)^m$ positive integers less than $N$ only divisible by those primes. Hence we see the image $\mathbb{N}$ under $q(x)$ is too dense to be only divisible by a finite set of primes.

If we apply our main theorem to such a sequence, choosing our isomorphism $\Pi_\U \overline{\mathbb{F}}_{p_n} \cong \mathbb{C}$ such that $lim_\U t_n = t$, then {\L}o\'s's theorem for ultraproducts (see \cite{Schoutens} section 2.3) tells us our algebra $A$ corresponds to a sequence of algebras $A_n$ in $Rep(S_{t_n})$ over $\overline{\mathbb{F}}_{p_n}$ for infinitely many values of $n$.

However, we know that that $A$ is inside the abelian subcategory of $\underline{Rep}(S_t)$ generated by the defining representation and its first $k$ tensor powers for some value of $k$ (this is true for all objects of $\underline{Rep}(S_t)$), and hence the same holds for the algebras $A_n$. Moreover the assumption that $A$ was not isomorphic to any of the algebras in the classification implies the same holds for $A_n$, which violates the above classification in $Rep(S_{t_n})$. $\square$

\medskip

\noindent \textbf{Remark:} One could similarly try to classify simple associative or Lie algebras in $\underline{Rep}(S_t)$, or more generally to classify simple $C$-algebras in $\underline{Rep}(S_t)$ for a linear operad $C$.  For a finite group $G$ Etingof gave a characterization of simple $G$-equivariant $C$-algebras in terms of subgroups of $G$ and simple (non-equivariant) $C$-algebras \cite{Etingof2}.  He then outlines how to use such a characterization over $\mathbb{C}$ to classify simple $C$-algebras in $\underline{Rep}(S_t)$ for transcendental $t$.  What comes out of our theory is that in order to carry out this program for algebraic $t$ one needs to have a characterization of simple $C$-algebras not just over $\mathbb{C}$, but also over algebraically closed fields in all sufficiently large characteristics.

\medskip

Going the other direction, we can take interesting observed phenomena in the setting of representation theory in complex rank to make new predictions about representation theory in large rank and characteristic.  We will close out this section by highlighting a question we think is particularly interesting in this direction:

\medskip

\noindent \textbf{Open problem:}  Entova-Aizenbud investigated representations of rational Cherednik algebras in complex rank $t$ in \cite{Inna}, where she found interesting degeneration phenomena at rational rank $t$.  What do these degeneration phenomena correspond to in the setting of rational Cherednik algebras in large characteristic and rank?

\end{subsection}

\begin{subsection}{Modular representations and supergroups}

Recently Entova-Aizenbud, Hinich, and Serganova constructed $\underline{Rep}(GL_n)^{ab}$ as a limit of the categories of representations of supergroups $GL(n+m|m)$ over $\mathbb{C}$ as $m$ tends to infinity \cite{EHS}.  On the other hand, our main theorem tells us that $\underline{Rep}(GL_n)^{ab}$ can be thought of as a limit of the categories $Rep(GL_{n+p}, \overline{\mathbb{F}}_p)$ as $p$ tends to infinity. 

We'll note that in the first case the limit is a filtered stable limit in the sense that there is an increasing filtration on $\underline{Rep}(GL_n)^{ab}$ such that each filtered piece is equivalent to a corresponding filtered piece of $Rep(GL(n+m|m))$ for all sufficiently large $m$, which is significantly stronger than the notion of limit involving an ultrafilter we have been using here. Nevertheless we can use these two interpretations of $\underline{Rep}(GL_n)^{ab}$ to pass results between $Rep(GL_{n+p}, \overline{\mathbb{F}}_p)$ for large $p$ and $Rep(GL(n+m|m))$ for large $m$.

In section 3.2 of \cite{EHS} they define a notion of admissible weights for $\mathfrak{gl}(n+m|m)$, corresponding to those highest weights of those representations occurring as subquotients of objects of the form $V^{\otimes p} \otimes V^{* \otimes q}$, where $V$ denotes the defining representation.  Combinatorially they are indexed by bipartitions $\lambda = (\lambda^{\circ}, \lambda^\bullet)$ where $\lambda^\circ$ has at most $n+m$ parts and $\lambda^\bullet$ has at most $m$ parts. One thing that comes out of their theory is that for a fixed bipartition $\lambda$, the ``admissible part" of the Kac module $K(\lambda)$ for $\mathfrak{gl}(n+m|m)$ stabilizes in an appropriate sense as $m \to \infty$, and corresponds to a standard object $V(\lambda)$ in $\underline{Rep}(GL_n)^{ab}$.

Highest weight theory for algebraic representations of $GL(n)$ is usually indexed by sequences $a_1 \ge a_2 \ge \dots \ge a_n$ of (possibly negative) integers.  However if we want to look at representations occurring as subquotients of objects of the form $V^{\otimes p} \otimes V^{* \otimes q}$ where $V$ is the defining representation and $n \gg p+q$ it is convenient to label them by bipartitions, where we encode a bipartition $\lambda = (\lambda^{\circ}, \lambda^\bullet)$ by the sequence $\lambda^{\circ}_1 \ge \lambda^{\circ}_2 \ge \dots \ge \lambda^\circ_\ell \ge0 \ge \dots \ge 0 \ge -\lambda^\bullet_m \ge -\lambda^\bullet_{m-1} \ge \dots \ge -\lambda^\bullet_1$.  In particular this reindexing corresponds to the usual labeling of irreducible objects by bipartitions in the Deligne category at generic values of $t$.

Now that we've fixed some indexing conventions, we are ready to state a proposition linking supergroup representations to modular representation theory:

\begin{proposition}
For fixed bipartitions $\lambda$, $\mu$ there exist $m_0$, $p_0$ such that whenever $m > m_0$ and $p > p_0$:

$$[K(\lambda),L(\mu)] = [\Delta(\lambda),L(\mu)],$$
where the left hand side is the multiplicity of the irreducible $GL(m+n | m)$ module $L(\mu)$ inside the Kac module $K(\lambda)$, and the right hand side is the multiplicity of the irreducible $GL(p+n, \overline{\mathbb{F}}_p)$ module $L(\mu)$ inside the Weyl module $\Delta(\lambda)$.
\end{proposition}

\noindent \textbf{Proof:} Both of these quantities coincide with the corresponding multiplicity of an irreducible inside a standard module inside the highest weight category $\underline{Rep}(GL_n)^{ab}$.  $\square$

 \end{subsection}

\end{section}

\end{document}